\documentclass[11pt]{amsart}

\numberwithin{equation}{subsection}
\usepackage{amsmath}
\usepackage{amsfonts}
\usepackage{amssymb}
\usepackage{amscd}
\newtheorem{conjecture} {Conjecture}
\newtheorem{theorem}{Theorem} [section]
\newtheorem{lemma} [theorem] {Lemma}
\newtheorem{proposition} [theorem] {Proposition}

\newtheorem{corollary} [theorem] {Corollary}
\newtheorem{definition} [theorem] {Definition}

\begin{document}
\title{Some $L^2$ results for $\overline\partial$ on projective varieties with general singularities}
\author{Nils \O vrelid and Sophia Vassiliadou }
\address{Dept. of Mathematics\\University of Oslo\\
P.B 1053 Blindern, Oslo, N-0316 NORWAY}
\address{Dept. of Mathematics\\Georgetown University\\
Washington, DC 20057 USA} \email{nilsov@math.uio.no,\;
sv46@georgetown.edu}
\maketitle
\date{\today}
\begin{abstract} \noindent Let $X$ be an irreducible $n$-dimensional
projective variety  in $\mathbb{CP}^N$ with arbitrary singular
locus. We prove that the
$L^2$-$\overline\partial$-$(p,1)$-cohomology groups (with respect
to the Fubini-Study metric) of the regular part of $X$ are finite
dimensional.

\end{abstract}

\medskip
\noindent

\section{Introduction}

\noindent This paper is the third in a series of papers in which
we discuss $L^2$ existence results for the Cauchy-Riemann operator
on  the regular part of  complex spaces with non-isolated
singularities. In \cite{FOV1} we considered the case of a
relatively compact, open, Stein subset $\Omega$ of a reduced Stein
complex space $X$ and studied the equation $\overline\partial u=f$
on $\Omega^*=\Omega\setminus A$ where $A$ was a lower dimensional
complex analytic set with empty interior containing the singular
locus of $X$ and $f$ was a $\;\overline\partial$-closed, square
integrable form vanishing to high order on $A$. In \cite{OV} we
proved some $L^2$ existence theorems for $\overline\partial$ in
the case of product singularities.

\medskip
\noindent In this paper we consider an irreducible $n$-dimensional
projective variety $X$ in $\mathbb{CP}^N$ with arbitrary singular
locus ${\mbox{Sing}}X$. The regular part of $X$, denoted by
${\mbox{Reg}}X$, inherits a metric from the restriction of the
Fubini-Study metric, which we call the ambient metric.  We denote
by $<,>_{FS},\; |\;|_{FS},\; dV_{FS}$ respectively the pointwise
inner product, norm on multi-(co) vectors and volume element on
${\mbox{Reg}}X$ induced by the ambient metric and by
$L^2_{p,q}({\mbox{Reg}}X, dV_{FS})$ the space of $(p,q)$ forms on
${\mbox{Reg}}X$ that are square-integrable with respect to the
ambient metric. We consider the weak $\overline\partial$-operator
acting on forms on ${\mbox{Reg}}\,X$. Let
$Z^{p,q}_{(2)}({\mbox{Reg}}X):=\{f\in L^2_{(p,q)}({\mbox{Reg}}X,
dV_{FS})|\;\; \overline\partial f=0 \;{\mbox{on}} \;
{\mbox{Reg}}X\}$ and let $f\in Z^{p,q}_{(2)}({\mbox{Reg}}X)$ where
$0\le p\le n$ and $1\le q\le n$. We address the question of
whether we can solve $\overline\partial u=f$ on ${\mbox{Reg}}X$
with $u\in L^2_{(p,q-1)}({\mbox{Reg}}X, dV_{FS} )$. Our main
result is the following theorem:

\begin{theorem}
There exists a closed subspace $\mathcal{H}$ of finite codimension
on the space $Z^{p,1}_{(2)}({\mbox{Reg}}X)$ \;such that for every
$f\in \mathcal{H}$ there exists a  $u\in L^2_{(p,0)}(Reg\,X,
dV_{FS})$ such that $\overline\partial u=f$\;on $RegX$.
\end{theorem}

\medskip
\noindent Pardon and Stern proved in \cite{PS1} that the
$L^2$-$(n,q)$-$\overline\partial$-cohomology groups of
$n$-dimensional projective varieties with general singularities
(with respect to the ambient metric) are finite dimensional for
all $q\ge 0$ (more precisely they showed that they are isomorphic
to the $(n,q)$-cohomology groups of a desingularization). To our
knowledge Theorem 1.1 is the first result on finite dimensionality
of some $L^2$-$(p,q)$-$\overline\partial$-cohomology groups (with
respect to the ambient metric) with $p\neq n$ of such varieties.
The case of projective varieties with isolated singularities has
been studied by Nagase \cite{N}, Ohsawa \cite{O1, O2, O3, O4, O5},
Pardon \cite{P} and Pardon and Stern \cite{PS1, PS2, PS3}. A local
analogue of this question for varieties with isolated
singularities was considered by Forn\ae ss \cite{F}, Diederich,
Fornaess and Vassiliadou \cite{DFV} and Forn\ae ss, \O vrelid and
Vassiliadou \cite{FOV1} (see \cite{DFV, FOV1} for extended
references). Ohsawa constructed in \cite{O6} a complete K\"ahler
metric on a small, Stein, deleted neighborhood of a singular point
of a complex space with arbitrary singularities and showed that
with respect to that metric, certain
$L^2$-$\overline\partial$-cohomology groups with support
conditions vanish there. Few years later, Grant and Milman
constructed in \cite{GM} a complete K\"ahler metric on the regular
part of a compact complex space with arbitrary singular locus.
Grant and Milman's construction was motivated in part by the hope
that this metric might be useful for proving the existence of a
Hodge structure on the intersection cohomology of projective
varieties with arbitrary singularities. It is not clear to us at
this moment whether the {\it{local}} $L^2$ results for
$\overline\partial$ obtained in \cite{O6} (with respect to the
specific {\it{complete}} K\"ahler metric) or Grant and Milman's
{\it{complete}}  K\"ahler metric can be used to understand the
{\it{global}} $L^2$-$\overline\partial$-cohomology groups (with
respect to the {\it{Fubini-Study}} metric) we consider in this
paper.

\medskip
\noindent Our proof is inspired by the methods employed in
\cite{DFV} and \cite{F}. It is based on the observation that
H\"ormander's $L^2$-theory carries over unchanged to Stein Riemann
domains $Y$ over $\mathbb{C}^n$ when we equip $Y$ with the pull
back Euclidean metric from $\mathbb{C}^n$. We look at the affine
pieces $X_0, X_1,\cdots, X_N$ of $X$ and for each piece we choose
a family of non-degenerate projections\;(for the definition see
Section 2) $\pi_j: X_i\to \mathbb{C}^n$ along with subvarieties
$\Sigma_{\pi_j}$ such that $\pi_j: X_i\setminus \Sigma_{\pi_j}\to
\mathbb{C}^n$ is a local biholomorphism and $\cup ( X_i\setminus
\Sigma_{\pi_j})={\mbox{Reg}}X$. Now the set  $Z_j:=X_i\setminus
\Sigma_{\pi_j}$ inherits a metric from the pull back of the
Euclidean metric in $\mathbb{C}^n$ and let us denote the pointwise
norm with respect to this metric by $|\;|_j$ and volume element
$dV_j$. Let $(\pi_j, Z_j, \Sigma_{\pi_j})_{j=1}^{M}$ be an
enumeration of the various projections of the affine pieces. When
$f\in L^2 ({\mbox{Reg}}X, dV_{FS})$ we can show that $f$ satisfies
estimates of the following form:

\begin{equation*}
\int_{Z_j} |f|^2_j \;e^{-\psi_j}\; dV_j \le C\;
\int_{{\mbox{Reg}}X} |f|_{FS}^2 \; dV_{FS},
\end{equation*}

\noindent where $\psi_j$ are specific plurisubharmonic functions
and $C$ some positive constant independent of $j$. Then, we may
use H\"ormander's $L^2$-theory to obtain $L^2$ solutions $v_j$
(with respect to the metric induced on $Z_j$ by the pull back of
the Euclidean metric in $\mathbb{C}^n$) to  $\overline\partial
v=f$ on $Z_j$. The delicate part is to use Lojasiewicz's
inequalities to obtain good control of the boundary behaviour of
$|\;|_j, \; dV_j$ on $Z_j$ in terms of the Fubini-Study metric.
The forms $h_{jj'}:=v_j-v_{j'}$ are holomorphic on $Z_j\cap
Z_{j'}$. From what we know about their boundary behaviour we can
deduce from sheaf-theoretic results that they lie in a finite
dimensional vector space of holomorphic $p$-forms. Since the map
$f\to \{h_{jj'}\}_{\{1\le j<j'\le M\}}$ is linear, we have that
$h_{jj'}=0$ for all $j,j'$, whenever  $f$ lies in a finite
codimensional subspace of $Z^{p,1}_{(2)}({\mbox{Reg}}X)$. In this
case, we can show that the $v_j$ define an $L^2$ solution $v$
(with respect to the ambient metric) to $\overline\partial v=f$ on
${\mbox{Reg}}X$.

\medskip
\noindent As a by-product of the techniques used in the paper we
obtain a weighted $L^2$-estimate for $\overline\partial$ on
irreducible affine algebraic subvarieties of $\mathbb{C}^N$. More
precisely we prove the following theorem:

\begin{theorem}Let $X$ be an irreducible, $n$-dimensional  affine subvariety of
$\mathbb{C}^N$ and let $\psi$ be a strictly plurisubharmonic
function on ${\mbox{Reg}}X$ with at most logarithmic growth (i.e.
$\psi(z)\le A\ log(1+\|z\|^2)+B$ for some $A,B\ge 0$) and not
necessarily bounded from below. Let $Z^{(p,1)}_{\psi}:=\{f\in
L^{2, loc}_{(p,1)}(RegX, dV_E);\;\; \overline\partial f=0\;\; on\;
RegX;\; \int_{RegX} |f|_E^2\;e^{-\psi} dV_E<\infty \}$. Then,
there exists a subspace $\mathcal{H}\subset Z^{(p,1)}_{\psi}$ of
finite codimension such that for all $f\in \mathcal{H}$ there
exists a $u\in L^{2,loc}_{(p,0)}(RegX, dV_E)$ with
$\overline\partial u=f$ on $RegX$ and

$$
\int_{Reg X} |u|^2_E \,(1+\|z\|^2)^{-2}\,e^{-\psi}\,dV_E\le C
\int_{RegX} |f|^2_E\,e^{-\psi}\,dV_E,
$$

\noindent where $C$ is some positive constant.
\end{theorem}

\medskip
\noindent Our main theorem can also be used to prove  finite
dimensionality results for some local
$L^2$-$\overline\partial$-cohomology groups of varieties with
isolated singularities. In particular, let $X$ be an irreducible
$n$-dimensional analytic set in $\mathbb{C}^N$ with an isolated
singularity at $0$. In \cite{FOV}\;(section 9)\, a question was
raised about understanding the
$L^2$-$\overline\partial$-$(p,q)$-cohomology groups (with respect
to the Euclidean metric) of the regular part of a small Stein
neighborhood of $0$ in $X$ when $p+q=n$ and $p,q>0$. With the aid
of Theorem 1.1 we can prove that when $p=n-1,\,q=1$ these groups
are finite dimensional.


\medskip
\noindent The paper is organized as follows: In section 2 we
describe some geometric facts about affine algebraic varieties and
recall some basic concepts about projective varieties. Section 3
deals with the comparison  between the Fubini-Study metric and the
pull-back metrics $|\;|_j$ that were defined earlier in the
introduction. Section 4 contains the estimates for the solutions
$v_j$ to $\overline\partial v=f$ on $Z_j$. In Section 5, we show
that $h_{jj'}$ lie in a finite dimensional subspace of the space
of holomorphic $(p,0)$ forms. In Section 6 we prove Theorem 1.2
and finally in section 7 we outline the proof for the finite
dimensionality of the local
$L^2$-$\overline\partial$-$(n-1,1)$-cohomology groups of a small
Stein neighborhood of an isolated singular point of an irreducible
$n$-dimensional variety of $\mathbb{C}^N$.

\medskip
\noindent {\bf{Acknowledgements:}} This paper was completed while
the second author was visiting the Institute of Mathematics at the
University of Oslo in June of 2006. She gratefully acknowledges
its hospitality and support. She would also like to thank Tom
Haines for many fruitful and stimulating discussions and Peter
Haskell for a valuable comment on an earlier version of the paper.

\section{Some geometric facts about  varieties}

\subsection{Non-degenerate projections}\;Let $Y$ be an irreducible, $n$-dimensional,  affine
algebraic variety in $\mathbb{C}^N$. Let $L$ be an $n$-dimensional
linear subspace of $\mathbb{C}^N$ and $\pi: \mathbb{C}^N \to L$ be
an orthogonal projection.

\begin{definition} We shall say that $\pi_{\upharpoonright Y}: Y\to L$
is a non-degenerate projection on $Y$, if there exists a point
$p\in {\mbox{Reg}}Y$ such that the induced map on tangent spaces
${(\pi_{\upharpoonright Y})}_{*, p}: T_p Y\to L$ is an
isomorphism.
\end{definition}

\medskip
\noindent{\bf{Remark:}}  Every projection $\pi:\mathbb{C}^N\to L,$
where $L$ is an $n$-dimensional subspace of $\mathbb{C}^N$ can be
made non-degenerate on $Y$ after a slight perturbation.

\begin{definition} The ramification locus $\Sigma$ of the above projection
consists of all points $p\in {\mbox{Sing}}Y$ as well as those
$p\in {\mbox{Reg}}Y$ for which ${(\pi_{\upharpoonright Y})}_{*,
p}: T_p Y\to L$ fails to be an isomorphism.
\end{definition}

\medskip
\noindent Another way to think about the non-degenerate
projections and their ramification loci is the following: The map
$\pi_{\upharpoonright Y}: Y\to L$ induces for all $y\in Y$, a
linear map $(\pi_{\upharpoonright Y})_{*,y}: {T_y} Y \to L$, where
${T_y}Y$ is the Zariski tangent space at $y$. As a linear map
between vector spaces it has a rank and a corank (the dimension of
the kernel of the linear mapping $(\pi_{\upharpoonright
Y})_{*,y}$). An equivalent way to restate  Definition 2.1 would be
that $\pi_{\upharpoonright Y}: Y\to L$ is non-degenerate on $Y$ if
there exists a regular point $p$ in $Y$ such that ${\mbox{rank}} (
(\pi_{\upharpoonright Y})_{*,p})=n$. Under this light we can
define the ramification locus of a projection
$\pi_{\upharpoonright Y}: Y\to L$ as follows:

\begin{definition} The ramification locus of the above projection
is the set

$$\Sigma=\{y\in Y;\;\;corank ((\pi_{\upharpoonright Y})_{*,y})\ge 1\;\}.
$$
\end{definition}

\medskip
\noindent

\medskip
\noindent \begin{proposition} The ramification locus $\Sigma$ is
an algebraic subvariety of $Y$. When the projection
$\pi_{\upharpoonright Y}: Y\to L$ is non-degenerate on $Y$ then
$\Sigma$ is a proper subvariety of $Y$.
\end{proposition}

\noindent \begin{proof}  The proof of the Proposition is similar
to that of Theorem 4 (page 136, Volume II) in Gunning \cite{G} (we
just replace holomorphic subvarieties by algebraic and holomorphic
map by regular map). The fact that $\Sigma$ is a proper subvariety
of $Y$ when the projection is non-degenerate on $Y$ follows
trivially from the definitions.\end{proof}

\medskip
\noindent{\bf{Remark:}} For a non-degenerate projection on $Y$,
the set $A:=\{y\in {\mbox{Reg}}Y;\;\;
{\mbox{rank}}((\pi_{\upharpoonright Y})_{*,y})=n\}$ is a non-empty
open Zariski dense set in $Y$. Thus, for most points $y\in
{\mbox{Reg}}Y$ the projection $\pi_{\upharpoonright Y}$ is
non-degenerate.

\subsection{Projections onto the $n$-coordinate planes}\;\;Let
$I:=(i_1,\cdots,i_n),\;1\le i_1<\cdots <i_n\le N$ be an increasing
$n$-tuple and let $I'$ be its complement in $\{1,2,\cdots,N\}$.
Let $L_I$ be an $n$-dimensional subspace of $\mathbb{C}^N$ defined
by $L_I:=\{(z_1,\cdots, z_N)\in
\mathbb{C}^N;\;\;z_j=0\;{\mbox{for\;all\;}}\; j\in I'\}$ and  such
that $\pi_I: \mathbb{C}^N \to L_I,\;\; (z_1,\cdots,z_N)\to
(z_{i_1},\cdots, z_{i_n}),$  the projection onto the $I$
coordinates is non-degenerate on $Y$ (we can always assume that
after a slight perturbation). Since $Y$ is an irreducible affine
variety in $\mathbb{C}^N$ there exists a prime ideal
$\mathcal{B}\subset \mathbb{C}[z_1,\,\cdots, z_N]$ such that
$Y=Z(\mathcal{B})$, the zero locus of $\mathcal{B}$. Let
$p_1,\cdots, p_s$ be generators of $\mathcal{B}$. It is a
well-known fact that the germs $p_{1,z}\cdots,p_{s,z}$ generate
$\mathcal{I}_{Y,z}$, the stalk at $z$ of the ideal sheaf of $Y$,
when the latter is viewed as a holomorphic variety; see for
example Proposition 13.3.3 in \cite{T}. Let $\Sigma_I$ denote the
ramification locus of $\pi_I$ and let $y\in \Sigma_I$. The induced
map on tangent spaces $({\pi_I}_{\upharpoonright Y})_{*,y}: T_y
Y\to L_I$ sends an element $T_y Y \ni v=(v_1,\cdots, v_N) \to
(v_{i_1},\cdots, v_{i_n})$. We know that the
${\mbox{kern}}(({\pi_I}_{\upharpoonright Y})_{*,y})$ is at least
one complex dimensional if and only if  all $(N-n)\times (N-n)$
minors of the matrix

$$\left(\frac{\partial p_i}{\partial z_r}(y)\right)_{\underset{ r\in I'}{1\le i\le
s,}}
$$

\noindent have zero determinants. Points $y\in {\mbox{Sing}}Y$ are
in $\Sigma_I$ since the rank of the Jacobian $\left(\frac{\partial
p_i}{\partial z_r}(y)\right)_{\underset{1\le r\le N}{ 1\le i\le
s,}}$ at these points is less that $N-n$.

\noindent
\medskip
\noindent

\begin{proposition} Let $Y$ be an irreducible,
$n$-dimensional, affine algebraic variety in $\mathbb{C}^N$. There
exist finitely many, complex linear, orthogonal projections
$\pi_{j}: \mathbb{C}^N \to L_j$, with $n$-dimensional images
$L_{j}$ and a constant $c>0$ such that:

\medskip
\noindent i)The restriction of $\pi_j$'s on $Y$  are
non-degenerate on $Y$ with ramification locus $\Sigma_j$.

\medskip
\noindent ii) For every $z\in {\mbox{Reg}}Y$ there exists a $j$
such that $z\notin \Sigma_j$ and such that $(\pi_j)_{*,z}: T_z
Y\to L_j$ satisfies $\|(\pi_j)_{*,z} v\|\ge c\;\|v\|$ for all
$v\in T_z Y$. (Here the norms are the induced Euclidean norms).

\medskip
\noindent In particular ii) implies that $\cap
\Sigma_j={\mbox{Sing}}Y$.
\end{proposition}

\begin{proof}  After a slight perturbation we can take as  $\pi_j$'s
the projections onto the (newly defined) $n$- dimensional
coordinate planes in $\mathbb{C}^N$. Let $I=\{i_1,\cdots,
i_n\},\;1\le i_1<\cdots<i_{n}\le N$ be an increasing $n$-tuple.
Let $\pi_{I}: C^N\to L_I$ be the projection onto the I
coordinates. To prove part ii) we need the following lemma.

\medskip
\noindent \begin{lemma} There exists an absolute  constant $c>0$
such that for each $z\in {\mbox{Reg}}Y$ there exists a multi-index
$I_0$ such that  $\|({\pi_{I_0}}_{\upharpoonright Y})_{*,z} v\|\ge
c \;\|v\|$ for all $v\in T_z Y$
\end{lemma}

\begin{proof} Let $z\in {\mbox{Reg}}Y$ and $v_1,\cdots v_n$ be an orthonormal basis of
$T_z Y$. Let $\{e^I:=e_{i_1}\wedge\cdots\wedge e_{i_n};|I|=n\}$ be
an orthonormal basis of $\bigwedge^n \mathbb{C}^N$. Consider the
expansion of $v_1\wedge v_2\cdots \wedge v_n$ in terms of $e^I$.
In what follows for abbreviation we shall write ${\pi_{I}}_*$
instead of $({\pi_I}_{\upharpoonright Y})_{*,z}$. We have

\begin{equation}
1=\|v_1\wedge\cdots\wedge v_n\|^2={\sum_{|I|=n}}'
\|{\pi_{I}}_{*}(v_1)\wedge \cdots\wedge {\pi_{I}}_{*}(v_n)\|^2=
{\sum_{|I|=n}}' D_I
\end{equation}

\noindent where
$D_I:=G({\pi_{I}}_{*}(v_1),\cdots,{\pi_{I}}_{*}(v_n))$ is the Gram
determinant of the vectors
${\pi_{I}}_{*}(v_1),\cdots,{\pi_{I}}_{*}(v_n)$.

\medskip
\noindent Recall that on a unitary space $E$ endowed with a
hermitian inner product $(\;,\;)$ the Gram determinant of vectors
$x_1,\cdots, x_p$ in $E$ is described by:

$$G(x_1,\cdots, x_p):={\mbox{det}} \begin{pmatrix} (x_1,x_1)\cdots &
(x_1,x_p)\\
\cdots & \cdots\\
(x_p,x_1)\cdots& (x_p, x_p)\\
\end{pmatrix}.$$

\medskip
\noindent In general $G(x_1,\cdots,x_p)\ge 0$ and equality holds
if the vectors $x_1,\cdots,x_p$ are linearly dependent.

\medskip
\noindent Since there are $\binom{N}{n}$ terms on the right hand
side of equation (2.2.1) there should exist a multi-index $I_0$
such that $D_{I_0} \ge {\binom{N}{n}}^{-1}$. Now if we let
$S:=\pi_{I_0}^*\circ \pi_{I_0}: T_z Y\to T_z Y$ we obtain

$$
D_{I_0}={\mbox{det}}\biggl((Sv_i,v_j)\biggr)={\mbox{det}}S.
$$

\noindent But $S$ is a positive, symmetric form that has
eigenvalues $\{l_j\},\; 0<l_1\le l_2\cdots \le l_n\le 1$ since
$\|S\|\le 1$. Then for all $v\in T_z Y$ we have

\begin{equation*}
\|\pi_{I_0} v\|^2=(Sv,v)\ge l_n \|v\|^2 \ge (\prod_{j} l_j)\;
\|v\|^2 =D_{I_0} \|v\|^2\ge {\binom{N}{n}}^{-1} \|v\|^2.
\end{equation*}
\end{proof}

\medskip
\noindent {\bf{Remarks:}}\;a)\;Part ii) of the above Proposition
can be thought of as  a statement about $n$-planes in
$\mathbb{C}^N$. Recall that the Grassmannian ${\mbox{Gr}}(n,N)$
can be covered by open affine sets $U_{\Gamma}$ where $\Gamma$ is
a $(N-n)$-dimensional subspace of $\mathbb{C}^N$. Each
$U_{\Gamma}$ is defined to be the subset of planes $\Lambda\subset
\mathbb{C}^N$ complementary to $\Gamma$. Fixing any subspace
$\Lambda\in U_{\Gamma}$ a subspace $\Lambda'\in U_{\Gamma}$ is the
graph of a homomorphism $\phi: \Lambda\to \Gamma$, so that
$U_{\Gamma}={\mbox{Hom}}(\Lambda,\Gamma)$\,(for more information
on this the interested reader may look at Lecture 16 in
\cite{Ha}). We shall consider a slightly different covering of the
Grassmanian ${\mbox{Gr}}(n,N)$. Let $L$ be an $n$-dimensional
subspace as before and let $L^{\perp}$ denote its orthogonal
complement in $\mathbb{C}^N$. Let $\mathcal{B}(L, {L}^{\perp})$
denote the set of bounded linear maps from $L$ to ${L}^{\perp}$.
Consider the following local parametrizations:

\begin{equation*}
\phi_{L}:\; \mathcal{B}(L,{L}^{\perp})\to Gr(n,N)
\end{equation*}

\noindent given by $\phi_L(T)={\mbox{Graph}}(T)$. Let
$\epsilon>0$. The set $\phi_{L}(\{T; \|T\|<\epsilon\})$ is an open
neighborhood of $L$ in ${\mbox{Gr}}(n,N)$. Since the latter
variety is compact there exist finitely many $L_1,\cdots,L_K\in
Gr(n,N)$ and open sets $U_j:=\phi_{L_j}(\{T; \|T\|<\epsilon\})$
such that $Gr(n,N)=\cup_{j=1}^K U_j$. When $L\in U_j,$ the
orthogonal projection $\pi_j: L\to L_j$ is bounded from below by
$(1+\epsilon^2)^{-\frac{1}{2}}$. Choosing $\epsilon$ small enough
we can make the constant $c$ that appears in part ii) of the above
proposition to be as close to 1 as we like, using sufficiently
many projections.
\end{proof}
\medskip
\noindent b) Part ii) of the above Proposition guarantees that
$\cap_{j} \Sigma_j ={\mbox{Sing}}Y$. Hence, $\{Y\setminus
\Sigma_j\}_{j}$ will cover ${\mbox{Reg}}Y$.

\medskip
\noindent
\subsection{Projective spaces and Fubini-Study metric} A point in $\mathbb{CP}^N$ is usually written as a homogeneous
vector $[Z_0,:\cdots,:Z_N]$ by which we mean the line spanned by
$(Z_0,\cdots,Z_N)\in \mathbb{C}^{N+1}\setminus\{0\}$. For
$i=0,\cdots,N$ we define

$$
\phi_i: \mathbb{C}^N\to U_i\subset \mathbb{CP}^N
$$

\medskip
\noindent given by $\phi_i(z_1,\cdots, z_N)=[z_1:\cdots :z_i:
1:z_{i+1}:\cdots:z_N]$.

\medskip
\noindent In particular $\phi_0(z_1,\cdots,z_N)
=[1:z_1:\cdots:z_N]$. We set for $0 \le i \le N$,
\;$$H_i:=\mathbb{CP}^N\setminus
\phi_i(\mathbb{C}^N)=\{[Z];\;Z_i=0\}.$$

\medskip
\noindent On $\phi_i^{-1}(U_i),$  using the affine coordinates
$z_1,\cdots,z_N$ the Fubini-Study metric takes the form

$$
\left(\sum h_{\mu\nu}(z) dz_{\mu}\otimes d\overline{z}_{\nu}
\right) (1+\|z\|^2)^{-2}
$$

\noindent where $h_{\mu\nu}(z)=(1+\|z\|^2)
\delta_{\mu\nu}-\overline{z}_{\mu}\,z_{\nu},\;\;\mu,\nu=1,\cdots,n$.
The associated $(1,1)$ K\"ahler form is described by

$$
\omega=\frac{i}{2} \sum h_{\mu\nu} dz_{\mu}\wedge
d\overline{z}_{\nu}.
$$

\medskip
\noindent Let $\hat{\lambda}_1,\cdots, \hat{\lambda}_N$ be the
eigenvalues of the restriction of the Fubini-Study metric on the
affine piece $\phi_i^{-1}(U_i)$ with respect to the Euclidean
metric. A direct calculation shows that

$$(1+\|z\|^2)^{-2}=\hat{\lambda}_1\le
\hat{\lambda}_2=\cdots=\hat{\lambda}_N=(1+\|z\|^2)^{-1}.$$

\section{Comparison
between various $L^2$-norms}

\medskip
\noindent Let $X$ be an irreducible $n$-dimensional projective
variety in $\mathbb{CP}^N$. Set $X_i:={\phi_i}^{-1}(X\cap
U_i)\subset \mathbb{C}^N$. For each $i,\,i\in \{0,\cdots,N\}$ we
shall choose $\{L^i_k\}_{k=1}^{M_i},$ families of $n$-dimensional
complex subspaces of $\mathbb{C}^N$ and orthogonal linear
projections $\pi^i_k: \mathbb{C}^N\to L^i_k$ such that
${\pi^i_k}_{\upharpoonright X_i}$ is non-degenerate on $X_i$ with
ramification locus $\Sigma^i_k$ and such that part ii) of
Proposition 2.4 holds for each $i$. Set $W^i_k:=X_i\setminus
\Sigma^i_k$. To reduce the number of indices we choose an ordering
of the set $\{(i,k); \; i\in \{0,\cdots, N\}, \;k=1,\cdots, M_i\}$
and of the corresponding objects $L^i_k,\;\pi^i_k,\;W^i_k$ such
that we have a bijection

\begin{eqnarray*}
\Theta: \{1,\cdots, M\}&\longrightarrow& \{(i,k)\in
\mathbb{N}^2;\;i\in\{0,\cdots,
N\}, 1\le k\le M_i\},\\
&j \longrightarrow& (i(j), k(j))\\
\end{eqnarray*}

\medskip
\noindent Let $L_1,\cdots, L_M$ be this ordering of all the
$n$-dimensional subspaces $\{L^i_k\}$ and let $\pi_1,\cdots,\pi_M$
be the corresponding projections. In what follows the index $i(j)$
will determine the affine variety $X_i$ that contains
$W_j,\;\Sigma_j$.

\medskip
\noindent Let $<,\,>,\;|\;|,\;dV$ denote the pointwise inner
product, norm on muti-(co)-vectors and volume element on
${\mbox{Reg}}{X}_i$ induced by the Fubini-Study metric, and let
$<\;>_E,\; |\;|_E,\;dV_E$ those induced on ${\mbox{Reg}}{X}_i$
from the Euclidean metric in $\mathbb{C}^N$ and $<,\,>_j,\;
|\;|_j,\; dV_j$ those pull-backed on $W_j\subset
{\mbox{Reg}}{X}_{i(j)}$ via $\pi_j$ from the Euclidean metric on
$L_j$. By the min-max principle we know that the eigenvalues of
the restriction of the Fubini-Study metric with respect to the
Euclidean metric on ${\mbox{Reg}}X_i$ satisfy

$$
(1+\|z\|^2)^{-2}\le \lambda_1\le
\lambda_2=\cdots=\lambda_n=(1+\|z\|^2)^{-1}. $$

\medskip
\noindent For $f\in L^{2,\,loc}_{p,q}({\mbox{Reg}}X_i, dV)$ we
have

\begin{equation}\label{eq:L2bound1}
(1+\|z\|^2)^{-\alpha}\;|f|^2_E\;dV_E\le |f|^2\;dV
\end{equation}

\noindent while for a form $u\in
L^{2,\,loc}_{p,q-1}({\mbox{Reg}}X_i,dV_E)$ we have

\begin{equation}\label{eq:L2bound2}
(1+\|z\|^2)^{-\beta} |u|^2\;dV\le |u|^2_E\;dV_E,
\end{equation}

\noindent where $\alpha,\,\beta$ are some non-negative constants
that depend only on $p,q$.

\subsection{Comparison between $dV_j, \;dV_E$ on $W_j$}

\medskip
\noindent The set $W_j$ inherits two metrics. One from the
restriction of the Euclidean metric in $\mathbb{C}^N$ and another
from the pull-back of the Euclidean metric on $L_j$ via the map
$\pi_j$. We begin this section by relating the volume elements of
these two metrics on $W_j$.

\begin{lemma} There exists a smooth function $m_j$ defined on $W_j$ such that

\medskip
\noindent i)  $dV_E=m_j\; dV_{j} \;\;{\mbox{on}}\;\; W_j.$

\medskip
\noindent ii) The function $\log m_j$ is plurisubharmonic on
$W_j$.
\end{lemma}

\begin{proof} Part i) can be taken as the definition of the function $m_j$. To prove part ii)
of the lemma we need a local description of the function $m_j$.
Without loss of generality we can  assume that the $n$-dimensional
subspace $L_j$ corresponds to the $n$-coordinate plane $L_I$ in
$\mathbb{C}^N$ where $I=(1,\cdots,n)$ is an increasing $n$-tuple.
Then $\pi_j:X_{i(j)}\to L_j$ is the projection onto the first
$n$-coordinates. The ramification locus of this projection is
characterized by the vanishing of the determinants of all
$(N-n)\times (N-n)$ minors of the matrix

$$\left(\frac{\partial p_i}{\partial z_j}\right)_{\underset{n+1\le j\le N}{1\le i\le s}}$$

\noindent  where  $p_1,\cdots, p_s$ are the generators of the
ideal of the variety $X_{i(j)}$. We can describe the set $W_j$ as

$$W_j= \{z_0\in X_{i(j)};\;\;\Delta_K(z_0)\neq 0 \;{\mbox{for\;some
\;multi-index K:}}\; 1\le k_1<\cdots<k_{N-n}\le s\},$$

\medskip
\noindent where $$\Delta_K(z_0):=\frac{\partial (p_{k_{1}},\cdots,
p_{k_{N-n}})}{\partial(z_{n+1},\cdots, z_N)}(z_0).$$

\medskip
\noindent Let $z_0\in W_j.$ Then there exists  a neighborhood $V$
of $z_0$ in $X_{i(j)}$ that is parametrized as
$(z',g_1(z'),\cdots,g_{N-n}(z'))$ for some functions $g_j$ and
with $z'=(z_1,\cdots,z_n)\in \pi_{i(j)}(V)$ (by the implicit
function theorem). For all $k\in K$ we have that
$p_k(z',g_1(z'),\cdots,g_{N-n}(z'))=0$. The implicit function
theorem allows us to compute for all $l$ with $1\le l\le N-n$ and
$\nu$ with $1\le \nu\le n$,

\begin{equation}\label{eq:pdofG}
\frac{\partial g_l}{\partial
z_{\nu}}(\pi_{i(j)}(z))=\frac{A_{\nu\widehat{K_l}}(z)}{\Delta_K(z)}
\end{equation}

\noindent where $A_{\nu\widehat{K_l}}$ is the determinant of the
$(N-n)\times (N-n)$ matrix $\left(\frac{\partial p_{k_i}}{\partial
z_{j}}\right)_{\underset{n+1\le j\le N}{1\le i\le N-n}}$ where the
$l$-th column has been replaced by $^{\top}(-\frac{\partial
p_{k_1}}{\partial z_{\nu}},\cdots,-\frac{\partial
p_{k_{N-n}}}{\partial z_{\nu}})$.

\medskip
\noindent For $i=1,\cdots,n$ and $z\in V$ we let

$$\zeta_i:=e_i+\sum_{k=1}^{N-n}\frac{\partial g_k}{\partial
z_i}(z') \;e_{n+k}
$$

\noindent where $\{e_i\}_{i=1}^N$ is the standard basis  of
$\mathbb{C}^N$. It is not hard to show that $\{\zeta_i\}_{i=1}^n$
form a basis of $T_{z} W_j$. Moreover for $1\le i\le n$ we have
${(\pi_{i(j)})}_{*,z} \zeta_i=\widetilde{e_i},$ where
$\widetilde{e_i}$ is the standard basis in $L_j$. Clearly

\begin{equation}\label{eq:volumecomp}
dV_E= {\mbox{det}} B\;dV_j,
\end{equation}

\medskip
\noindent where $B=(b_{kl})$ is the $n\times n$ matrix with
entries $b_{kl}:=<\zeta_{k},\,\zeta_l>_E$  and  $<,\,>_E$ is the
pointwise Euclidean inner product on elements of $T_z W_j$.

\medskip
\noindent Let us look at $\|\bigwedge_{i=1}^n \zeta_i\|^2_E$. It
follows from the definition of the Euclidean inner product on
vectors in $\bigwedge^{n} T_z W_j$ that $\|\bigwedge_{i=1}^n
\zeta_i\|^2_E={\mbox{det}}B$. Using $(\ref{eq:volumecomp})$ we
have the following local description of
$m_j(z)=\|\bigwedge_{i=1}^n \zeta_i\|^2_E$.

\medskip
\noindent Let us also consider the following $n\times N$ matrix:

$$
C=
\begin{bmatrix}
1&     \cdots&    0&    \frac{\partial g_1}{\partial z_1}&     \cdots&      \frac{\partial g_{N-n}}{\partial z_1}\\
\cdot&  \cdots&   \cdot&         \cdot&                                 \cdots& \cdot\\
0&        1&      0&     \frac{\partial g_1}{\partial z_i}&     \cdots&      \frac{\partial g_{N-n}}{\partial z_i}\\
\cdot&   \cdots&   \cdot&     \cdot&               \cdots& \cdot\\
0&     \cdots&    1&          \frac{\partial g_1}{\partial z_n}&     \cdots&      \frac{\partial g_{N-n}}{\partial z_n}\\
\end{bmatrix}
$$

\medskip
\noindent To prove ii) we notice that $m_j(z)= {\sum}'_{H}
|C_H|^2,$ where the summation runs over all multi-indices
$H=(h_1,\cdots,h_n)$ with $1\le h_1<\cdots<h_n\le N$ and $C_H$ are
determinants from the $n\times N$ matrix $C$ with columns
$h_1,\cdots,h_n$. The $C_H$'s are holomorphic functions. Hence
$\log m_j$ is a plurisubharmonic function in a neighborhood of
$z_0$.
\end{proof}

\medskip
\noindent We have seen that $\{\zeta_j\}_{j=1}^n$ form an
orthonormal basis for $T_z W_j$ with respect to the pull back
metric. Let $\{\lambda_j\}$ be the eigenvalues of the K\"ahler
form $\omega_E$ of the Euclidean metric on $W_j$ with respect to
the pull-back metric. We have $1\le \lambda_1\le \cdots\le
\lambda_n$. For a $(p,q)$ form $f$ on $W_j$ we have:

\begin{equation}\label{eq:defnptnorm}
|f|^2_E={\sum_{I,J}}' \frac{1}{\prod_{i\in I}
\lambda_i\;\;\prod_{j\in J} \lambda_j} |f_{IJ}|^2.
\end{equation}

\medskip
\noindent Taking into account part i) of the previous lemma and
the fact that the eigenvalues $\{\lambda_j\}$ are greater or equal
to 1 we have the following estimate for $(p,q)$- forms $f$ with
$p>0,\;q>0$

\begin{equation}\label{eq:compjande}
m_j^{-1}\;|f|^2_j dV_j\le |f|^2_E\;dV_E
\end{equation}

\medskip
\noindent In the special case where $p=0,\;q>0$ we have a stronger
estimate

$$
|f|^2_j\;dV_j\le |f|^2_E\;dV_E.
$$

\subsection{Global Lojasiewicz inequality} In the next section we
shall need an upper bound for the function $m_j$ that was defined
in the previous lemma. The upper bound will be obtained by
applying a global Lojasiewicz-type inequality obtained first by
Brownawell \cite{B} and later improved by Ji, Koll\'ar and
Shiffman \cite{JKS}.

\begin{corollary} (Corollary 6 in \cite{JKS}) \noindent Let
$f_1,\cdots,f_k\in \mathbb{C}[z_1,\cdots, z_n]$ and let $d_i=
{\mbox{deg}}f_i$. Let $Z\subset \mathbb{C}^n$ be the common zero
set of these polynomials. Then there is a constant $C>0$ such that

\begin{equation}\label{eq:GLI}
\left(\frac{{\mbox{dist}}(z,
Z)}{(1+\|z\|^2)}\right)^{\overline{B}(n,\; d_1,\cdots,d_k)}\le C\;
{\mbox{max}}_i\;|f_i(z)|
\end{equation}

\noindent holds for all $z\in \mathbb{C}^n$. Here $C$ and
$\overline{B}(n, d_1,\cdots,d_k)$ are positive constants that
depend on the $f_i$ and $\|z\|^2:=|z_1|^2+\cdots+|z_n|^2$.
\end{corollary}

\medskip
\noindent  Recall that $m_j(z)={\mbox{det}}\left(
<\zeta_k,\,\zeta_l>_E\right)_{k,l}$ and
$<\zeta_k,\,\zeta_l>_E=\delta_{kl}+ \sum_{\mu=1}^{N-n}
\frac{\partial g_{\mu}}{\partial
z_k}\overline{\left(\frac{\partial g_{\mu}}{\partial
z_l}\right)}$. Using $(\ref{eq:pdofG})$ we can obtain the
following upper bound for $m_j(z)$:

\begin{equation}\label{eq:roughbound}
m_j(z)\le \left(1+|\Delta_K|^{-2} \sum_{\nu=1}^{n}\sum_{l=1}^{N-n}
|A_{\nu\widehat{K_l}}|^2\right)^n
\end{equation}

\medskip
\noindent The ramification locus $\Sigma_j$ is the zero locus of
the following set of polynomials $\{p_1,\cdots, p_s\}\cup
\{\Delta_K; {\mbox{for\;all\; increasing\; multi-indices\;}} K,\;
|K|=N-n\}$. Applying $(\ref{eq:GLI})$ we have for $z\in W_j$

\begin{equation}\label{eq:boundforD_K}
{\mbox{max}}_K |\Delta_K (z)|\ge {C'}_{i(j)}^{-1}\; d_E(z,
\Sigma_{j})^{\overline{B}_j}\;(1+\|z\|^2)^{-\overline{B}_j}
\end{equation}

\medskip
\noindent If we choose for each $z\in W_{j}$ the multi-index $K$
such that $\Delta_K(z)$ is maximal then using
$(\ref{eq:boundforD_K})$ inequality $(\ref{eq:roughbound})$ will
become

\begin{equation}\label{eq:betterbound}
m_j(z)\le \left(1+
{C'}_{i(j)}^2\;d_E(z,\,\Sigma_{j})^{-2\overline{B}_j}\;(1+\|z\|^2)^{2\overline{B}_j}\;
\;\sum_{\nu=1}^{n}\sum_{l=1}^{N-n}
|A_{\nu\widehat{K_l}}|^2\right)^n
\end{equation}

\medskip
\noindent Finally there exist constants $C_{i(j)},\; E_{i(j)},\;
D_{j}>0$ such that

\begin{equation}\label{eq:finalboundmj}
m_j(z)\le
C_{i(j)}\;(1+\|z\|^2)^{D_j}\;\left({\mbox{min}}\{1,\;d_E(z,\Sigma_{j})\}\right)^{-E_{i(j)}}
\end{equation}




\section{$L^2$-solvability for $\overline\partial$ on  $W_{j}$}

\medskip
\subsection{$L^2$-existence theorem}
\noindent H\"ormander's  $L^2$ theory for Stein domains in
$\mathbb{C}^n$ extends naturally to Stein Riemann domains over
$\mathbb{C}^n$ when these are given the pull back metric. In our
case $(W_j, \pi_j, L_j)$ are Riemann domains over $L_j$. Since the
sets $X_i$ are singular though, the $W_j$'s are not necessarily
Stein. We shall need the following variant of H\"ormander's
theory:

\begin{proposition} Let $q>0$ and  $f\in L^{2,loc}_{p,q}(W_j, dV_E)$ with
$\overline\partial f=0$ on $W_j$. Then there exists a solution $u$
to $\overline\partial u=f$ on $W_j$ satisfying the following
estimate:

\begin{equation}\label{eq:L2estSR}
\int_{W_j} |u|^2_j\; (1+\|\pi_j(z)\|^2)^{-2}\;e^{-\psi}\;dV_j\le
\int_{W_j} |f|^2_j \;e^{-\psi}\;dV_j,
\end{equation}

\noindent whenever the RHS is finite and $\psi$ is a
plurisubharmonic function on $W_j$.
\end{proposition}

\begin{proof} We can choose polynomials $Q_j(z)$ such that they vanish on
$\Sigma_j$ but do not vanish identically on $X_{i(j)}$. Let
$Z(Q_j)$ denote the zero set of $Q_j$. Then $(X_{i(j)}\setminus
Z(Q_j), \pi_j, L_j)$ is a Stein Riemann domain. Theorem 4.4.2 in
\cite{H} (or Theorems  2.2.1 and $2.2.1^{\prime}$ in \cite{H1})
carries (resp. carry) over almost verbatim to
$W'_j:=X_{i(j)}\setminus Z(Q_j)$ and we obtain the existence of a
solution $u\in L^{2,\ loc}_{p,q-1}(W'_j, dV_j)$ to
$\overline\partial u=f$ on $W'_j$ satisfying (\ref{eq:L2estSR})
with $W_j$ being replaced by $W'_j$. We want to show that the
solution $u$ extends to $W_j$ and satisfies a similar estimate
there. Let us look at the set $W_j\setminus W'_j=W_j\cap Z(Q_j)$.
It is a hypersurface in the complex manifold $W_j$. If we could
show that $u$ is locally in $L^2$ (with respect to $dV_j$) near
$W_j\setminus W'_j$ then the result would follow from the
following lemma:

\begin{lemma} Let $\Omega$ be an open subset of $\mathbb{C}^n$
and $Y$ an analytic subset of $\Omega$. Assume that $v$ is a
$(p,q-1)$ form with $L^{2,loc}$ coefficients and $w$ a
$(p,q)$-form with $L^{1,loc}$ coefficients such that
$\overline\partial v=w$ on $\Omega\setminus Y$ (in the sense of
distributions). Then $\overline\partial v=w$ on $\Omega$.
\end{lemma}

\begin{proof}  This is Lemma 6.9 (page 485) in \cite{D}.
\end{proof}

\medskip
\noindent {\bf{Remark:}} To be more precise we need to prove a
similar  lemma where $\Omega$ is an open set of a complex
manifold. But, Lemma 4.2 generalizes quite easily in this case.

\medskip
\noindent Now for every $z_0\in W_j\setminus W'_j$ there exists a
neighborhood $U$ of $z_0$ in $X_{i(j)}$ such that $\pi_j: U\to
\pi_j(U)\subset \mathbb{C}^n$ is a biholomorhism and there exists
a positive constant $c_1$ such that
$e^{-\psi(z)}\,(1+\|\pi_j(z)\|^2)^{-2}\ge c_1$ for all $z\in U$.
Hence $u\in L^2(U\setminus Z(Q_j), dV_j)$ and $\overline\partial
u=f$ on $U\setminus Z(Q_j)$. But then $\overline\partial u=f$ on
$U$ and hence $\overline\partial u=f$ on $W_j$.
\end{proof}

\medskip
\noindent \subsection{Estimates for the solution to
$\overline\partial u=f$ on $W_j$}  In what follows, we shall use
the notation $L^{2}_{p,q}(W, \phi, dV_h)$ to denote the Hilbert
space of $(p,q)$ forms on a complex hermitian manifold $(W,h)$ for
which $\|f\|^2:=\int_W |f|^2_{h} \; e^{-\phi} dV_{h}<\infty$. We
shall also use the notation $a\lesssim b$ (resp. $a\gtrsim b$) if
there exists an absolute positive constant $c$ such that $a\le c\;
b$ \;(resp. $a\ge c \;b$),\; $a\approx b$ if there exist absolute
positive constants $c, c'$ such that $a\le c \;b,\; b\le c'\; a$.

\medskip
\noindent We shall choose as $\psi_j(z):=\log m_j(z)+\alpha\,
\log(1+\|z\|^2)$ where $\alpha$ is the positive constant that
appears in (\ref{eq:L2bound1}). Then
$e^{-\psi_j}=m_j^{-1}\;(1+\|z\|^2)^{-\alpha}$. Let $f\in
L^2_{p,q}(\ Reg X, dV_{FS})$ with $\overline\partial f=0$ there.
Then ${\phi^*_{i(j)}}f\in L^2_{p,q}(\ Reg X_{i(j)}, dV)$. Using
(\ref{eq:L2bound1}),\;(\ref{eq:compjande})  we obtain

$$
m_j^{-1}\;(1+\|z\|^2)^{-\alpha}\;|{\phi^*_{i(j)}} f|^2_j\;dV_j\le
|{\phi^*_{i(j)}}f|^2\;dV.
$$

\medskip
\noindent But this last inequality implies that
${\phi^*_{i(j)}}f\in L^2_{p,q}(W_j, \psi_j, dV_j)$. Then we can
apply Proposition 4.1 to obtain a solution $u'_j\in
L^{2,loc}_{p,q-1}(W_j)$ that satisfies

$$
\int_{W_j} |u'_j|^2_j\;
(1+\|\pi_j(z)\|^2)^{-2}\;e^{-\psi_j}\;dV_j\le \int_{W_j}
|{\phi^*_{i(j)}}f|^2_j \;e^{-\psi_j}\;dV_j,
$$

\medskip
\noindent Taking into account (\ref{eq:L2bound2}), Lemma 3.1 i)
and the fact that $\|\pi_j(z)\|\le \|z\|$ we obtain

\begin{equation}\label{eq:prebound}
\int_{W_j} m_j^{-2}\;(1+\|z\|^2)^{-2-\alpha-\beta}\;|u'_j|^2\;dV
\le \int_{\ Reg X} |f|_{FS}^2 dV_{FS}
\end{equation}

\medskip
\noindent Using (\ref{eq:finalboundmj}) the last inequality
becomes

\begin{equation}\label{eq:almostfinbound}
\int_{W_j} (1+\|z\|^2)^{-2-\alpha-\beta-D_j}\; ({\mbox{ min}}\{1,
d_E(z,\Sigma_j)\})^{E_{i(j)}}\; |u'_j|^2\; dV \le \int_{\ RegX}
|f|_{FS}^2\;dV_{FS}.
\end{equation}

\medskip
\noindent We set $\Sigma_j^*:=\phi_{i(j)} (\Sigma_j)\cup (X\cap
H_{i(j)})$. Then ${\Sigma^*_{j}}$ is a projective subvariety of
$X$. Recall that $d_E(z,
\Sigma_j)>d(\phi_{i(j)}(z),{\Sigma^*_{j}})$ and $d(\phi_{i(j)}(z),
H_{i(j)})\approx (1+\|z\|^2)^{-\frac{1}{2}}$ where by $d(\bullet,
\bullet)$ we denote the projective distance. Taking into account
these inequalities we obtain

\begin{equation}\label{eq:almfinalbound}
\int_{{\phi_{i(j)}}^{-1}(\phi_{i(j)}(W_j))} d(\phi_{i(j)}(z),
H_{i(j)})^{D'_j}\;({\mbox{min}}\{1,
d(\phi_{i(j)}(z),\;\Sigma_j^*)\})^{E_{i(j)}}\;|u'_j(z)|^2\;dV(z)\le
\int_{Reg X} |f|^2_{FS}\;dV_{FS}
\end{equation}

\noindent where $D'_{j}:=2\,(2+\alpha+\beta+D_j)$. To relate
$d(\phi_{i(j)}(z), H_{i(j)})$ to $d(\phi_{i(j)}(z), \Sigma^*_j)$
we need to recall the notion of regular separation. Following
Lojasiewicz (page 242 in \cite{L}) we define:

\begin{definition} Let $(E,F)$ be closed subsets of a manifold M.
We say that $(E,F)$ satisfy the condition of regular separation if
 for each point $w\in E\cap F$ the following inequality holds
true in a neighborhood of the point $w$

$$
\rho(z,E)+\rho(z,F)\ge c \;\rho(z, E\cap F)^p
$$

\noindent where  $c,p$ are some positive constants.
\end{definition}

\medskip
\noindent
\begin{theorem} Every pair of analytic sets of a complex
manifold satisfies the condition of regular separation.
\end{theorem}

\begin{proof} This theorem is proved in \cite{L} (page 244).
\end{proof}

\medskip
\noindent Since the sets $X,\;H_i$ are regularly separated we have
for $z\in X_i$

$$
d(\phi_{i(j)}(z),H_{i(j)})\ge c\, d(\phi_{i(j)}(z), H_{i(j)}\cap
X)^p\;\ge c\,d(\phi_{i(j)}(z),\Sigma_j^*)^{p_{i(j)}}.
$$

\noindent where $c,p$ are the positive constants that appear in
the definition of regular separation.

\bigskip
\noindent Pulling back to $X$ via the biholomorphism
${\phi_i}^{-1}: U_i\to \mathbb{C}^N$  we obtain from
$(\ref{eq:almfinalbound})$

\begin{equation}\label{eq:almostfinalbound1}
\int_{\phi_{i(j)}(W_j)} ({\mbox{min}}\{1,d
(z,\Sigma^*_j)\})^{b_j}\; | u_j|_{FS}^2\;dV_{FS} \le \int_{\ RegX}
|f|_{FS}^2\;dV_{FS},
\end{equation}

\noindent where $b_{j}:=p_{i(j)}\,D'_{j}+E_{i(j)}$ and
$u_j:={({\phi_{i(j)}}^{-1})}^{*} u'_j$. Now, since $\mathbb{CP}^N$
has a finite diameter with respect to the Fubini-Study metric we
have that ${\mbox{min}}\{1, d(z, \Sigma^*_j)\}\approx
d(z,\Sigma_j^*)$. Hence $\overline\partial u_j=f$ on
$\phi_{i(j)}(W_j)$ and from $(\ref{eq:almostfinalbound1})$ we see
that  $u_j$ satisfies

\begin{equation}\label{eq:finalbound}
\int_{\phi_{i(j)}(W_j)} d (z,\Sigma^*_j)^{b_j}\; |
u_j|_{FS}^2\;dV_{FS} \le \int_{\ RegX} |f|_{FS}^2\;dV_{FS}.
\end{equation}

\medskip
\noindent To summarize: We have found solutions $u_j$ that satisfy
$\overline\partial u_j=f$ in $\phi_{i(j)} (W_j)$ and the estimate
$(\ref{eq:finalbound})$. For each $j,j'$ with $1\le j<j'\le M$ we
have $\overline\partial \left( u_j-u_{j'}\right)=0$ on
$\phi_{i(j)}(W_j)\cap \phi_{i(j')}(W_{j'})$. Now, the sets
$\phi_{i(j)}(W_j)$ are nonempty, Zariski open in $X$ whose
complement in $X$ is $\Sigma^*_j$. Set $\Sigma:=\cup_{j=1}^M
\Sigma^*_j$. We can restrict $u_j$ on $X\setminus \Sigma$. Let
$h_{jj'}:=(u_j-u_{j'})_{\upharpoonright (X\setminus \Sigma)}$.
Using $(\ref{eq:finalbound})$ we can show that the $h_{jj'}$'s
satisfy the following estimates:

\begin{equation}\label{eq:finalbound1}
\int_{X\setminus \Sigma} d (w, \Sigma)^{b'_{jj'}}
|h_{jj'}|^2_{FS}\,dV_{FS} \lesssim \int_{\ RegX}
|f|_{FS}^2\;dV_{FS}<\infty
\end{equation}

\noindent where $b'_{jj'}$ are some positive constants that depend
on $b_{j},\,b_{j'}$.

\subsection{Construction of the global solution}
\noindent Let us define the linear mapping

\begin{eqnarray*}
T: Z^{p,1}_{(2)}({\mbox{Reg}} X)&\to& {(\Omega^p(X\setminus
\Sigma))}^{\binom{M}{2}}\\
 f&\to& (h_{jj'})_{1\le j<j'\le M}\\
\end{eqnarray*}
\smallskip
\noindent Suppose we could show that the range of $T$ were finite
dimensional. Then  $\mathcal{H}:={\mbox{kern}}T$ would be a
subspace of $Z^{(p,1)}_{(2)}({\mbox{Reg}}X)$ of finite
codimension. For an $f\in \mathcal{H}$ we shall have that
$u_j=u_{j'}$ on $X\setminus \Sigma$ which is non-empty, Zariski
open subset of $\phi_{i(j)}(W_j)\cap \phi_{i(j')}(W_{j'})$. But
then $u_j=u_{j'}$ on $\phi_{i(j)}(W_j)\cap \phi_{i(j')}(W_{j'})$.
To finish the proof of the main theorem it would suffice to show
that the $u_j$'s determine a global, square-integrable (with
respect to the ambient metric) $(p,0)$-form on ${\mbox{Reg}}X$.

\medskip
\noindent Let $R$ be a positive real number with $R>\sqrt{N+1}$.
Then $\mathbb{CP}^N=\cup_{i=0}^N \phi_i (B(0,R)),$ where $B(0,R)$
is the Euclidean ball centered at the origin in $\mathbb{C}^N$ and
having radius $R$. Consider the sets $ F_j:=\{z\in
X_{i(j)}\setminus \Sigma_j; \;\|z\|<R\;{\mbox{and}}\;
\|({\pi_j})_{*,z} v\|\ge c \|v\|\;{\mbox{for\;all}}\;v\in T_z
X_{i(j)}\}.$ Let $V_j:=\phi_{i(j)}(F_j)\subset \phi_{i(j)}(W_j)$.

\begin{lemma} The sets $V_j$ cover ${\mbox{Reg}}X$.
\end{lemma}

\begin{proof} Obvious from the definition of the sets $V_j$.
\end{proof}

\noindent When $z\in F_j$ we have the following upper bound for
$m_j(z)$; $m_j(z)\le c^{-2n}$. Hence by $(\ref{eq:prebound})$ we
obtain that the solutions $u'_j$ satisfy the following
$L^2$-estimate on $F_j$'s:

$$
\int_{F_j} c^{-2n}\;(1+R^2)^{-2-\alpha-\beta} |u'_j|^2\;dV\le
\int_{{\mbox{Reg}}X} |f|^2_{FS}\;dV_{FS}
$$

\smallskip
\noindent Pulling back to $X$ via $\phi_{i(j)}^{-1}$ we see that
the solutions $u_j$ satisfy an $L^2$-estimate (with respect to the
ambient metric) on $V_j$. Since $u_j=u_{j'}$ on $V_j\cap V_{j'}$
and $\{V_j\}_j$ cover  ${\mbox{Reg}}X$ we obtain a global $(p,0)$
form $u$ that satisfies $\overline\partial u=f$ on ${\mbox{Reg}}X$
and

$$
\int_{{\mbox{Reg}}X} |u|_{FS}^2\;dV_{FS}\le C\int_{{\mbox{Reg}}X }
|f|^2_{FS}\;dV_{FS}.
$$

\medskip
\noindent The next section is devoted to proving that the
${\mbox{Rang}}(T)$ is finite dimensional.

\section{Finite dimensionality of certain holomorphic
$p$-forms}

\medskip
\noindent The main goal in this section is to prove the following
lemma:

\begin{lemma} Let $A$ be a nonnegative real number.  Let

$$
\mathcal{E}_A:=\{h\in \Omega^p(X\setminus \Sigma);\;\;\int_{w\in
X\setminus \Sigma} d^A(w,{\Sigma})\; |h(w)|_{FS}^2\;
dV_{FS}(w)<\infty\}.
$$

\noindent For any $A\ge 0, \;\mathcal{E}_A$ is a finite
dimensional complex vector space.
\end{lemma}

\noindent
\begin{proof} The idea of the proof is to ``kill'' the singularity of $h\in
\mathcal{E}_A$ by tensoring it with a section of an invertible
sheaf on $X$ vanishing to high order on $\Sigma$ and then use
classical finiteness results for $\Gamma(X,\mathcal{S}_k)$ where
$\mathcal{S}_k$ is a suitable coherent analytic
$\mathcal{O}_X$-module.

\medskip
\noindent Let $p:\mathbb{C}^{N+1}\setminus \{0\}\to \mathbb{CP}^N$
be the standard projection map. In what follows we shall think
mostly of $X$ as a compact complex space. The sheaves on $X$ that
we will consider shall be analytic sheaves. By $\mathcal{O}_X$ we
shall denote the sheaf of holomorphic functions on $X$.  For every
positive integer $l$ we define the twisting sheaf $\mathcal{O}(l)$
on $\mathbb{CP}^N$ as follows: if $U$ is an open subset of
$\mathbb{CP}^N$ in the Euclidean topology then $\mathcal{O}(l)(U)$
consists of the space of holomorphic functions on $p^{-1}(U)$
which are homogeneous of degree $l$. The global holomorphic
sections of this sheaf can be naturally identified with
homogeneous polynomials $Q(Z_0,\cdots, Z_N)$ of degree $l$ in
$\mathbb{C}^{N+1}$. Over $U_i,$ $\mathcal{O}(l)$ has a
trivialization with transition functions
$\left(\frac{Z_j}{Z_i}\right)^{l}$ on $U_i\cap U_j$.

\medskip
\noindent We choose a homogeneous polynomial $P(Z)$ of degree $d$
that vanishes on $\Sigma$ but does not vanish identically on $X$.
Let $k\in \mathbb{N}$ and let $\sigma^k$ be the section of
$\mathcal{O}(kd)$ that corresponds to $P(Z)^k$. We shall prove
that when $h\in \mathcal{E}_A$ and $k$ is a sufficiently large
positive integer, $h\otimes \sigma^k$ extends to a global section
$\widetilde{h\otimes \sigma^k}$ of some suitably chosen coherent
analytic $\mathcal{O}_X$-module $\mathcal{S}_k$. More precisely we
shall choose as $\mathcal{S}_k:=R^0 \pi_*
(\Omega^p_{\tilde{X}})\bigotimes_{\mathcal{O}_X}
\mathcal{O}(kd)_{\upharpoonright X}$\;  where $\pi:
\widetilde{X}\to X$ is a desingularization of $X$ such that
$\widetilde{\Sigma}:=\pi^{-1}(\Sigma)$ is a divisor with normal
crossings, \;$\Omega^p_{\tilde{X}}$ is the sheaf of holomorphic
$p$-forms on $\tilde{X}$ and $\mathcal{O}(kd)_{\upharpoonright
X}=i^{-1}{({\mathcal{O}}{(kd)})}$  with  $i: X\hookrightarrow
\mathbb{CP}^N$.

\medskip
\noindent Let us cover $X$ by finitely many open affine balls
$\mathcal{U}_{\nu}$ with $\mathcal{U}_{\nu}\subset\subset X\cap
U_{i}$ where $i=i(\nu)\in \{0,\cdots,N\}$. We shall first prove
that $h\otimes
{\sigma^k}_{\upharpoonright{\mathcal{U}_{\nu}\setminus \Sigma}}$
extends to a section $\mathcal{S}_k(\mathcal{U}_{\nu})$. Let us
work with inhomogeneous coordinates. Then the section $\sigma^k$
is represented by $p_i(z)=Z_i^{-kd} \;P(Z)^k$ for $z\in U_i$ and
we have $|p_i(z)^k|\le C \;d_E(z,\Sigma)^k$ for $z\in
\mathcal{U}_{\nu}$. When $h\in \mathcal{E}_A,$ $h\otimes \sigma^k$
is represented over $\mathcal{U}_{\nu}\setminus \Sigma$ by
$h\,p_i^k$ and

\begin{equation}\label{eq:1estimate}
\int_{z\in \mathcal{U}_{\nu}\setminus \Sigma} |h(z)\;p_i(z)^k|^2\;
d_E(z, \Sigma)^{A-2k}\;dV(z)<\infty
\end{equation}

\noindent The following pointwise estimates were proven in Lemma
3.1 in \cite{FOV1}:

\begin{lemma} We have for $x\in \pi^{-1}(\mathcal{U}_{\nu}\setminus \Sigma),\;
v\in \wedge^{r}T_x (\pi^{-1}(\mathcal{U}_{\nu}))$

\begin{eqnarray} \label{eq:pointwise}
c'\; {d}^t (x,\tilde{\Sigma})&\le& d(\pi(x), \Sigma)\le C'\; d(x,\tilde{\Sigma}),\\
c\; d^M (x, \tilde{\Sigma}) \;|v|_{x,\sigma} &\le&
|\pi_{*}(v)|_{\pi(x)}\le C\; |v|_{x,\sigma}.
\end{eqnarray}\\
\noindent  for some  positive constants $c',c,C',C,t,M,$ where
$c,C,M$ may depend on $r$ and where $\sigma$ is a real analytic
hermitian metric on $\pi^{-1}(X)$.

\medskip
\noindent For an $r$-form $a$ in $\mathcal{U}_{\nu}\setminus
\Sigma$ set\; $|\pi^{*}a|_{x,\sigma}:= {\mbox{max}} \{\;|<
a_{\pi(x)}, \pi_{*}v>|\;;\; |v|_{x,\,\sigma} \le 1,\; v\in
\wedge^{r}T_{x}(\pi^{-1}(\mathcal{U}_{\nu})\setminus
\widetilde{\Sigma})\}$, where by $<,>$ we denote the pairing of an
$r$-form with a corresponding multi-vector. Using (5.0.3) we
obtain:

\begin{equation}\label{eq:pointwiseest}
c\; d^{M}(x,\tilde{\Sigma}) \ \; |a|_{\pi(x)}\le |\pi^{*}
a|_{x,\sigma}\le C\;|a|_{\pi(x)}
\end{equation}

\noindent on $\pi^{-1}(\mathcal{U}_{\nu})$, for some positive
constant $M$.
\end{lemma}

\noindent Using the above lemma and choosing $k$ to be
sufficiently large we can show that $\pi^*(h\;p_i^k)\in
L^2_{(p,0)}(\pi^{-1}(\mathcal{U}_{\nu}\setminus \Sigma))$ and thus
extends to a holomorphic $p$-form on
$\pi^{-1}(\mathcal{U}_{\nu})$. But then
$h\otimes\sigma^k_{\upharpoonright{\mathcal{U}_{\nu}\setminus
\Sigma}}$ extends uniquely to a section in
$\mathcal{S}_k(\mathcal{U}_{\nu})$. The local extensions fit
together to a unique global section $\widetilde{h\otimes
\sigma^k}\in \Gamma(X,\mathcal{S}_k)$ extending
$h\otimes\sigma^k$. Since the map $h\to \widetilde{(h\otimes
\sigma^k)}$ is injective and $\Gamma(X, \mathcal{S}_k)$ is finite
dimensional, the space $\mathcal{E}_A$ is finite dimensional as
well.

\medskip
\noindent {\bf{Remark:}} H. Flenner showed in \cite{Fl} that
holomorphic $p$-forms on the regular part of a projective variety
with general singularities are the push-forward of holomorphic
$p$-forms on $\tilde{X}$, a desingularization of $X,$ provided
that $0\le p<{\mbox{codim Sing}} X-1$. In our case we do not need
to put any restriction on the range of $p$ since we are tensoring
our section $h$ with sections of an invertible sheaf on $X$.
\end{proof}

\subsection{Conclusion of proof of Theorem 1.} Clearly we can find an $A>0$
such that for all $1\le j<j'\le M,$ \;$h_{jj'}\in \mathcal{E}_A$.
Combining this with lemma 5.1 we have shown that the
${\mbox{Rang}}(T)$ is finite dimensional. Then the main theorem
follows by the argument used in Section 4.3.

\section{Proof of Theorem 1.2} \noindent One interesting application of the
techniques used in this paper is a derivation of a weighted
$L^2$-estimate for solutions to $\overline\partial u=f$ on the
regular part of affine varieties. Our goal in this section is to
prove Theorem 1.2.




\begin{proof} According to Proposition 2.4 there exist $\{L_j\}_{j=1}^M$ a
finite collection of $n$-dimensional subspaces of $\mathbb{C}^N$
and orthogonal projections $\pi_j: \mathbb{C}^N\to L_j$ such that:

\medskip
\noindent i) ${\pi_j}_{\upharpoonright(X)}$ is non-degenerate on
$X$ with ramification locus $\Sigma_j$ and

\medskip
\noindent  ii) For every $z\in {\mbox{Reg}}Y$ there exists a $j$
such that $z\notin \Sigma_j$ and such that $(\pi_j)_{*,z}: T_z
Y\to L_j$ satisfies $\|(\pi_j)_{*,z} v\|\ge c\;\|v\|$ for all
$v\in T_z Y$. (Here the norms are the induced Euclidean norms).

\medskip
\noindent Let us denote by $W_j:=X\setminus \Sigma_j$. Then
$(W_j,\pi_j, L_j)$ are Riemann domains and we would like to apply
H\"ormander's $L^2$ theory for $\overline\partial$ (in particular
Proposition 4.1) to our $f_{\upharpoonright(W_j)}$. Let us define
on $W_j$ the function $\psi_j:=\psi+{\mbox{log}} m_j$. Then
$\psi_j$ is plurisubharmonic on $W_j$ and using (3.1.4) we can
show that

$$
\int_{W_j} m_j^{-1} |f|_j^2 \,e^{-\psi}\;dV_j\le \int_{W_j}
|f|_E^2\,e^{-\psi}\,dV_E<\infty
$$

\noindent Noticing that $m_j^{-1}\,e^{-\psi}=e^{-\psi_j}$ we
conclude that $f\in L^2_{(p,1)}(W_j, \psi_j, dV_j)$. Hence we can
apply Proposition 4.1 to each $W_j$ and obtain a solution $v_j$ to
$\overline\partial v=f$ on $W_j$ and $v_j$ satisfies the following
estimate:

\begin{equation}\label{eq:estforvj}
\int_{W_j} |v_j|_j^2\,
(1+\|\pi_j(z)\|^2)^{-2}\,e^{-\psi_j}\,dV_j\le C\int_{W_j}
|f|^2_E\,e^{-\psi}\,dV_E<\infty
\end{equation}

\noindent Noticing that $\|\pi_j(z)\|\le \|z\|$ we derive the
following estimate from $(\ref{eq:estforvj})$

$$
\int_{W_j} |v_j|^2_j\,(1+\|z\|^2)^{-2}\,e^{-\psi_j}\,dV_j\le
C\int_{W_j} |f|^2_E\,e^{-\psi}\,dV_E<\infty
$$

\medskip
\noindent We want to estimate $v_j$ using the restriction of the
Fubini-Study metric on $Reg X$. Recall that $|v_j|_j\ge |v_j|_E,$
$dV_j=m_j^{-1}dV_E$ and $m_j\le C_{j}
(1+\|z\|^2)^{D_j}\,{\mbox{min}}\{1,d_E(z, \Sigma_j)\}^{-E_j}$.
Hence $|v_j|^2_j\, dV_j\ge C^{-1}_{j}
(1+\|z\|^2)^{-D_j}\,{\mbox{min}}\{1,\,d_E(z,\Sigma_j)\}^{E_j}\,|v_j|^2_E\,dV_E$.
Taking into account that $|v_j|^2_E\,dV_E\ge
(1+\|z\|^2)^{-\beta}\,|v_j|^2\,dV$ we obtain

\begin{equation}\label{eq:FSestvj}
e^{-\phi_j}\,|v_j|^2_j \,dV_j\ge e^{-\psi}C^{-2}_{j}
(1+\|z\|^2)^{-2\,D_{j}-\beta}\,{\mbox{min}}\{1,\,d_E(z,\Sigma_j)\}^{2\,E_j}\,|v_j|^2\,dV.
\end{equation}

\medskip
\noindent We can consider the projective closure $PCl(X)$ of $X$
in $\mathbb{CP}^N$ and let $\phi_0: \mathbb{C}^N\to
U_0\subset\mathbb{CP}^N$ be the corresponding affine chart. Let
$\Sigma^*_j:=\phi_0(\Sigma_j)\cup (PCl(X)\cap H_0)$ and let
$\Sigma:=\cup \Sigma^*_j$. Recall that
$(1+\|z\|^2)^{-\frac{1}{2}}\,\approx d(\phi_0(z), H_0)\ge
C'\,d(\phi_0(z),\Sigma^*_j)^{p}$ and
$d_E(z,\Sigma_j)>d(\phi_0(z),\Sigma^*_j)$. As in section 4 we can
show that

$$
(1+\|z\|^2)^{-2\,
D_j-\beta}\,{\mbox{min}}\{1,\,d_E(z,\Sigma_j)\}^{2\,E_j}\ge C''
d^{2(2\,D_j+\beta)\,p+2\,E_j}(\phi_0(z), \Sigma).
$$

\medskip
\noindent Hence the $v_j$ satisfy the following estimate on $W_j$

\begin{equation}\label{eq:almfinestvj}
\int_{W_j} d(\phi_0(z), \Sigma)^{N'_j}\,e^{-\psi}\, |v_j|^2\,
dV\le C\int_{W_j} |f|^2_E\,e^{-\psi}\,dV
\end{equation}

\noindent where $N'_j$ are some positive constants. Since $\psi$
has at most logarithmic growth we can bound from below
$e^{-\psi}\ge e^{-B}\,(1+\|z\|^2)^{-A}\gtrsim d(\phi_0(z),
H_0)^{2A}\gtrsim d(\phi_0(z), \Sigma^*_j)^{2Ap}$. Hence we have

\begin{equation}\label{eq:finestvj}
\int_{W_j} d(\phi_0(z), \Sigma)^{N_j}\, |v_j|^2\, dV\le
C\int_{W_j} |f|^2_E\,e^{-\psi}\,dV
\end{equation}

\noindent where $N_j$ some positive constants that depend on
$N'_j, A$.

\medskip
\noindent Let us look at $u_j:=(\phi_0^{-1})^* v_j$. Then,
$\overline\partial (u_j-u_k)=0$ on $\phi_0(W_j)\cap \phi_0(W_k)$.
Now, each of the sets $\phi_0(W_j)$ is a non-empty, Zariski open
in $PCl(X)$ whose complement is $\Sigma^*_j$. Let us denote by
$h_{jk}:=(u_j-u_k)_{\upharpoonright(PCl(X)\setminus \Sigma)}$.
Then $\{h_{jk}\}_{1\le j<k\le M}$ are holomorphic $p$-forms on
$PCl(X)\setminus \Sigma$ and from $(\ref{eq:finestvj})$ we can see
that they satisfy the following estimates

\begin{equation}\label{estforh}
\int_{w\in PCl(X)\setminus \Sigma} d(w,
\Sigma)^K\;|h_{jk}(w)|^2_{FS}\,dV_{FS}(w)<\infty,
\end{equation}

\noindent where $K$ is a positive constant that depends on $N_j,
N_k$. But then, lemma 5.1 tells us that the space $\mathcal{E}_K$
of such forms is a finite dimensional complex vector space. Hence
we can repeat the argument in section 4.3 and construct a linear
operator $T: Z^{(p,1)}_{\psi}\ni f \to (h_{jk})_{1\le j<k\le M}\in
\mathcal{E}_K^{\binom{M}{2}}$. Then $\mathcal{H}:={\mbox{kern}}T$
is a finite codimensional subspace of $Z^{(p,1)}_{\psi}$. For an
$f\in Z^{(p,1)}_{\psi}$, the local solutions $\{u_j\}$ agree on
$PCl(X)\setminus \Sigma$, thus they agree on
$\phi_0(W_j)\cap\phi_0(W_k)$. Hence $v_j=v_k$ on $W_j\cap W_k$ and
we can define a $(p,0)$-form $u$ on $\cup W_j={\mbox{Reg}}X$ that
would satisfy $\overline\partial u=f$ on ${\mbox{Reg}}X$.
Repeating a similar argument as in the end of section 4.3 we can
show that this $u$ satisfies the desired estimate.
\end{proof}

\medskip
\noindent {\bf{Remarks:}} i) The codimension of $\mathcal{H}$ in
$Z^{(p,1)}_\psi$ will increase with $A$ the constant that appears
in the logarithmic growth for $\psi$ (and in general it will not
be finite). Indeed, let us consider the variety $X:=Y\times
\mathbb{C}$ where $Y$ is an irreducible surface in $\mathbb{C}^3$
described by a homogeneous polynomial of degree $d\ge 3$. In
\cite{F}, Forn\ae ss constructed finitely many
$\overline\partial$-closed $(0,1)$ forms $\lambda$ on the
${\mbox{Reg}}Y$  satisfying: 1) $\lambda\in L^2_{0,1}(Reg Y\cap
B_1(0), dV_E)$ and 2) $\overline\partial v=\lambda$ is not
solvable in $L^2(Reg Y\cap B_1(0))$ (here $B_1(0)$ is the unit
ball in $\mathbb{C}^N$ centered around $0$). For one of these
forms $\lambda$ one can further show that

\begin{equation}\label{eq:estforl}
\int_{z\in Reg Y} |\lambda(z)|^2_E\,
(1+\|z\|^2)^{-a}\,dV_E(z)<\infty
\end{equation}

\noindent when $a>1$.

\medskip
\noindent We also have

$$
\int_{w\in \mathbb{C}} |w|^{2k}\,(1+|w|^2)^{-b}\,dA(w)<\infty
$$

\noindent when $b>k+1$.

\medskip
\noindent Let $m\in \mathbb{N}$. On ${\mbox{Reg}}X$ we consider
forms $\lambda(z)\otimes p_m(w)$ where $p_m(w)$ is a polynomial in
$w$ of degree $m$ and $\lambda$ is the non-solvable $(0,1)$- form
on ${\mbox{Reg}}Y$ as above that satisfies $(\ref{eq:estforl})$.
Let us define the function $\psi:=A\,\log (1+\|(z,w)\|^2)$ where
$A:=a+\beta$ and where $a$ is as above and $\beta$ is chosen such
that $\beta>m+1$ (hence $A>m+2$). Clearly $ \lambda(z)\otimes
p_m(w)\in Z^{(0,1)}_{\psi}\ominus\mathcal{H}$ \,(otherwise the
equation $\overline\partial v=\lambda$ would have an $L^2$
solution in a deleted neighborhood of $0$ in $Y$ which would
contradict the choice of $\lambda$). Hence the codimension of
$\mathcal{H}$ in $Z^{0,1}_{\psi}$ will be greater or equal to
$m+1$ and it will be infinite whenever  $\psi$ has faster growth
than a logarithmic one.

\medskip
\noindent ii) The following question was posed to the first author
by Henkin and Zeriahi: Could it be possible to take as
$\mathcal{H}$ in Theorem 1.2 the whole space $Z^{(p,1)}_{\psi}$
when $X$ is non-singular?

\section{Another application of Theorem 1.1}

\medskip
\noindent Let $X$ be an irreducible $n$-dimensional variety in
$\mathbb{C}^N$ with an isolated singularity at $0$. Let $\Omega$
be a small Stein neighborhood of $0$ with $\partial\Omega$ smooth.
In this section we shall prove that the
$L^2$-$\overline\partial$-$(n-1,1)$-cohomology group of
${\mbox{Reg}}\,\Omega$ (with respect to the Euclidean metric)  is
finite dimensional. In what follows by
$L^2_{p,q}({\mbox{Reg}}\,\Omega)$ we denote the space of $(p,q)$
forms on ${\mbox{Reg}}\,\Omega$ that are square-integrable with
respect to the restriction of the Euclidean metric on
${\mbox{Reg}}\,\Omega$.

\medskip
\noindent We shall need the following general result:

\begin{proposition} Under the above assumptions and for any $p,q\in
\mathbb{N}$ with $q>0$ there exists a finite codimensional
subspace $E_0$ of $L^2_{(p,q)}(Reg\,\Omega)\cap
kern(\overline\partial)$ and a linear operator $S:E_0\ni f\to
u_0\in L^{2,\,loc}_{p,q-1}(\overline{\Omega}\setminus \{0\})$ such
that $\overline\partial u_0=f$ on ${\mbox{Reg}}\,\Omega$.
\end{proposition}

\medskip
\noindent Suppose for the moment that Proposition 7.1 were true.
We choose $\chi\in C^{\infty}(\overline{\Omega})$ such that
$\chi=1$ near $\partial\Omega$ and $\chi=0$ near $0$. For an $f\in
E_0$ we write $f=\overline\partial
(\chi\,u_0)+\overline\partial((1-\chi)\,u_0)$. Let
$f_0:=\overline\partial
((1-\chi)u_0)=(1-\chi)f-\overline\partial\chi\wedge u_0$. Then
$f_0\in L^2_{p,q-1}(Reg\Omega)$ and is $\overline\partial$-closed
there. Now $\Omega$ can be embedded as a subdomain of an
irreducible $n$-dimensional projective variety $X'$ and we may
extend $f_0$ by zero to $\overset{o}{f_0}$, a globally defined
form on ${\mbox{Reg}}X'$  that is $\overline\partial$-closed on
${\mbox{Reg}}X'$ and square-integrable (with respect to the
Fubini-Study metric) there. Let us consider the case where $q=1$.
According to Theorem 1.1 there exists a subspace $\mathcal{H}$ of
$Z^{p,1}_{(2)}(Reg X')$ of finite codimension such that whenever
$g\in \mathcal{H}$ we can solve $\overline\partial v=g$ with $L^2$
estimates (with respect to the Fubini-Study metric) on
${\mbox{Reg}}X'$. We apply Theorem 1.1 to $\overset{o}{f_0}$ and
we obtain a solution $v$ to $\overline\partial v=\overset{o}{f_0}$
on ${\mbox{Reg}}X'$.

\medskip
\noindent Let us consider the map $T: E_0\to Z^{p,1}_{(2)}(Reg
X')$ sending an element $ E_0\ni f\to \overset{o}{f_0}$. Clearly
$T$ is a linear map, hence ${\mbox{codim}}_{E_0}
(T^{-1}(\mathcal{H}))\le
{\mbox{codim}}_{Z^{p,1}_{(2)}(RegX')}(\mathcal{H})$. Setting
$E:=T^{-1}(\mathcal{H})$ we see that $E$ is of finite codimension
in $L^2_{p,1}(Reg\,\Omega)\cap {\mbox{kern}}(\overline\partial)$
and when $f\in E$ we can find a
$u:=\chi\,u_0+v_{\upharpoonright(Reg\,\Omega)}\in
L^2_{p,0}(Reg\,\Omega)$ satisfying $\overline\partial u=f$ on
${\mbox{Reg}}\,\Omega$. Taking $p=n-1$ we obtain the finite
dimensionality of the
$L^2$-$\overline\partial$-$(n-1,1)$-cohomology group of
${\mbox{Reg}}\,\Omega$.

\medskip
\noindent We return now to the proof of Proposition 7.1.

\begin{proof} We consider a desingularization $\pi: \tilde{X}\to
X$ with exceptional divisor $D$ and let $\sigma$ be a hermitian
metric on $\widetilde{X}$ and $dV_{\sigma}$ the volume element
induced by this metric. Let $\mathcal{O}(kD)$ denote the
holomorphic line bundle on $\tilde{X}$ associated to the divisor
$kD$ and let $\widetilde{\Omega}:=\pi^{-1}(\Omega)$. Choose  a
hermitian metric $h$ on $\mathcal{O}(kD)$. When $f\in
L^2_{p,q}(Reg\, \Omega),$\; $\pi^*f$ does not necessarily belong
to $L^2_{p,q}(\widetilde{\Omega},\,dV_{\sigma})$. However it gives
rise to a section $\xi_f\in L^2_{p,q}(\widetilde{\Omega},
 \mathcal{O}(kD)_{\upharpoonright {\tilde{\Omega}}})$ for some sufficiently large integer
 $k$. Moreover, $\overline\partial \xi_f=0$ on $\widetilde{\Omega}$.
 Since $\widetilde{\Omega}$ has strictly pseudoconvex boundary
 the $L^2$-$\overline\partial$-cohomology groups $H^{p,q}_{(2)}(\widetilde{\Omega},\,\mathcal{O}(kD))$
for $q>0$ are finite dimensional (see for example Theorem 5.11 in
\cite{KR}). Hence there
 exists a finite codimensional subspace $E$ of
 $L^2_{p,q}(\widetilde{\Omega},\,\mathcal{O}(kD))\cap
 {\mbox{kern}}(\overline\partial)$ such that whenever $\xi_f\in E$
 there exists $\tau_f\in L^2_{p,q-1}(\widetilde{\Omega},\,
 \mathcal{O}(kD))$ such that $\overline\partial \tau_f=\xi_f$.
 Choosing  $\tau_f$ to be the minimal solution we make the map
 $\xi_f\to \tau_f$ linear. Now ${\tau_f}_{\upharpoonright {\widetilde{\Omega}\setminus D}}$
 determines a form $\widetilde{v}_f\in
 L^{2,loc}_{p,q-1}(\overline{\widetilde{\Omega}}\setminus\{D\})$. Setting
 $u_0:=(\pi^{-1})^{*}\widetilde{v}_f$ we obtain a $u_0\in
 L^{2,loc}_{p,q-1}(\overline{\Omega}\setminus\{0\})$ satisfying
 $\overline\partial u_0=f$ on ${\mbox{Reg}}\,\Omega$. Moreover the
 set $E_0:=\{f\in L^{2}_{p,q}({\mbox{Reg}}\,\Omega)\cap
 {\mbox{kern}}(\overline\partial): \xi_f\in E\;\}$ is a finite
 codimensional subspace of $L^{2}_{p,q}({\mbox{Reg}}\,\Omega)\cap
 {\mbox{kern}}(\overline\partial),$ since the map $f\to \xi_f$ is
 injective.
\end{proof}

\medskip
\noindent {\bf{Remark:}} Arguments similar to the one used in the
proof of  Proposition 7.1 have already appeared in the works of
Nagase \cite{N}, Pardon \cite{P} and Pardon and Stern \cite{PS3}
(sections 3,4 of their 1997 preprint).

\end{document}